\documentclass{llncs}
\usepackage{graphicx}
\usepackage{epsfig}

\def\RR{{\rm I\hspace{-0.50ex}R}}

\title{Evolutionary Mesh Numbering: \\
Preliminary Results
}
\author{Francis Sourd  and Marc Schoenauer}
\institute{
CMAP -- URA CNRS 756\\
Ecole Polytechnique\\
Palaiseau 91128, France\\
{\tt marc.schoenauer@polytechnique.fr}\\~\\
in I. Parmee, ed., Proceedings of ACDM'98, pp 137--150, Springer Verlag, 1998
}
\date{}

\begin{document}

\maketitle

\begin{abstract}
Mesh numbering is a critical issue in 
Finite Element Methods, as the computational cost of one analysis is
highly dependent on the order of the nodes of the mesh.
This paper presents some preliminary investigations on the problem of
mesh numbering using Evolutionary Algorithms. 
Three conclusions can be drawn from these experiments. First, the
results of the up-to-date 
method used in all FEM softwares (Gibb's method) can be consistently
improved; second, none of the crossover operators tried so far (either
general or problem specific) proved useful; third, though the general
tendency in Evolutionary Computation seems to be the hybridization
with other methods (deterministic or heuristic), none of the presented
attempt did encounter any success yet.
The good news, however, is that this algorithm allows an
improvement over the standard heuristic method between 12\% and 20\% 
for both the 1545 and 5453-nodes meshes used as test-bed.
Finally, some strange interaction between the selection scheme and the
use of problem specific mutation operator was observed, which appeals
for further investigation.
\end{abstract}

\section{Introduction}
Most Design Problems in engineering make an intensive use of numerical
simulations of some physical process in order to predict the actual
behavior of the target part. When the mathematical model supporting
the numerical simulation involves Partial Differential Equations, 
Finite Element Methods (FEM) are today one of the most widely used
method by engineers to actually obtain an approximate numerical
solution to their theoretical model.
However, the computational cost of up-to-date
numerical methods used in FEM directly depends on 
the way the nodes of the underlying mesh (i.e. the discretization) are
numbered.  
Solving the Mesh Numbering Problem (MNP) amounts to find the
permutation of the order of the nodes that minimizes that computational
cost. 

Numerical engineers have
developed a powerful heuristic technique (the so-called {\em Gibb's
method}) that gives reasonably good
results, thus providing a clear
reference to any a new mesh numbering algorithm.

The goal of this paper is to use Evolutionary Algorithms (EAs) to
tackle the MNP. EAs have been widely
applied to other combinatorial optimization problems, among which the popular
TSP \cite{Grefenstette_TSP87,TSP:PPSN96}. However, as far as
we know, this is the first attempt to solve the MNP using EAs.
Unfortunately, though both problem look for a solution in the space of
permutations of $[0,n]$, the specificity of
the MNP might make inefficient the simple transposition of TSP
techniques to the MNP.
Indeed, looking at the history of Evolutionary TSP, it seems clear that
the key of success is hybridization with problem-specific techniques:
from the  Grefenstette's early incorporation of domain knowledge 
\cite{Grefenstette_TSP87} to the most recent works
 \cite{TSP:PPSN96,Merz:ICEC97,Asparagos:ICEC97} where evolutionary
results can -- at last -- be compared to 
the best Operational Research results, even for large size TSP 
instances. 
So the  path to follow here seems rather clear: 
design some NMP-specific operators, and compare their performances to
either ``blind'' problem independent operators or TSP-specific
operators.

The paper is organized as follows: Section \ref{FEM} recalls the basics
of Finite Element Methods, and precisely defines the objective
function. The state-of-the-art ``Gibbs method' is also briefly described.
Section \ref{EA} presents the design of the particular Evolutionary
Algorithm used thereafter, discussing in turn the representation
issue, specific crossover and mutation operators, 
and the initialization procedure.
This algorithm is then experimented in section \ref{results}, with
emphasis on the tuning of the 
probabilities of application of all operators at hand (both
problem-specific and problem independent). First, the crossover
operators all rapidly appear harmful. Second,   
surprising interactions between the selection scheme and the different
mutation operators seem to indicate that, though at the moment domain
knowledge did not increase the performances of the algorithm, there is
still some room for improvement in that direction.
Finally, the usefulness of evolutionary mesh numbering  with respect to 
 Gibbs' method is discussed: the results are indeed better, but the
cost is also several orders of magnitude greater. 

\section{Mesh Numbering}
\label{FEM}

\subsection{Theoretical Background}
Many models for physical, mechanical, chemical and biological 
phenomenon end up in Partial Differential Equations (PDE)
where the unknown are functions defined on some domain $\Omega$ of $\RR^n$. 

\begin{figure}
\begin{center}
\begin{tabular}{cc}
\psfig{figure=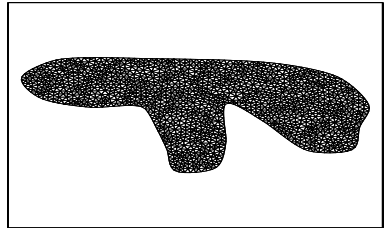,height=3cm,width=5cm}
\end{tabular}
\caption{\em Sample mesh with 1545 nodes.}
\label{test-mesh}
\end{center}
\end{figure}

A popular way to transform a system of PDEs into a finite linear
system is the Finite Element Method (FEM) \cite{Zienkiewicz:FEM,Ciarlet}. 
The domain $\Omega$ is discretized into small {\em elements}, who build up
a {\em mesh}. A typical mesh -- on a non-typical domain -- is given in
Figure \ref{test-mesh}. The solution of the initial
PDEs is sought in spaces of functions that are polynomial
on each element. Such an approximate
solution is completely determined by its values at some
points of each element, called the {\em nodes} of the mesh.
Those values are computed by writing the original PDEs locally on
each element, resulting in a linear system of size the number of nodes
times the number
of {\em degrees of freedom}, or unknown values, at each node). 
Usual sizes for such systems range from a few hundreds (e.g.
in two-dimensional structural mechanics) to millions
(e.g. for three-dimensional problems in aerodynamics).

However, as a consequence of 
the local discretization, the equation
at each node only involves values at a few neighboring nodes: the resulting
matrix is hence very sparse. And specific 
methods exist for sparse system \cite{Golub}, whose 
complexity is proportional to the square of the 
{\em bandwidth} of the matrix, i.e.
the average size of the {\em profile} of the matrix, given by the
sum over all lines of the maximal distance to the diagonal of non-zero
terms. For full matrices, the bandwidth is $n(n-1)/2$ while it is $n$
for tridiagonal matrices (for $n \times n$ matrices).

\vspace{-1cm}

\subsection{Computing the bandwidth}
\label{example}
The contribution of each single line to the total bandwidth
of the matrix is highly dependent on the order of the nodes
of the mesh: the equation for node number $i$ only involves 
the neighboring nodes; hence the only non-zero terms
of the corresponding equation will appear in the matrix  
in the column equal to the number of the node in the mesh.
Depending on the order of the nodes in the mesh, the bandwidth
can range from a few units to almost the size of the matrix.

\begin{figure}
\begin{center}
\begin{tabular}{cc}
\psfig{figure=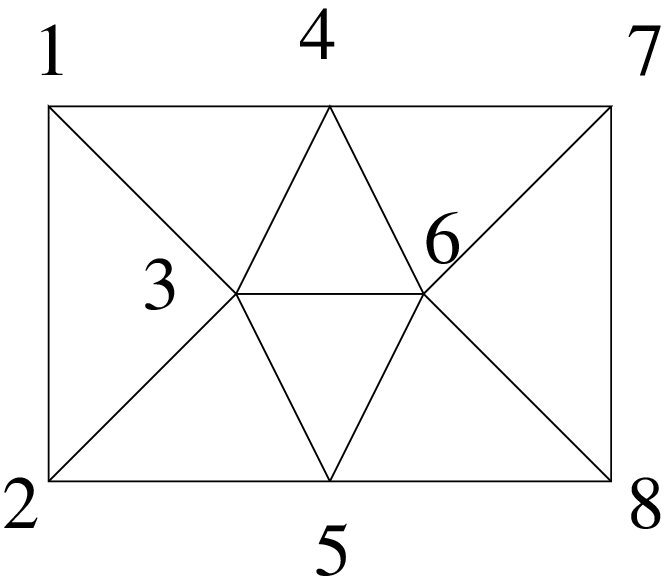,width=3cm} &
\psfig{figure=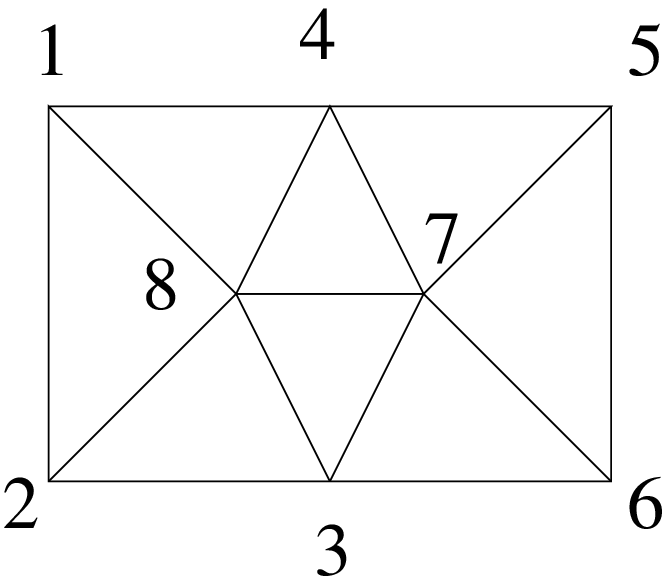,width=3cm} \\
(a) {\em Bandwidth=18} &
(b) {\em bandwidth=26}
\end{tabular}
\caption{\em A simple example of bandwidth with respect to nodes order.}
\label{ref_mesh}
\end{center}
\end{figure}

A simple example of a mesh is
given in Figure \ref{ref_mesh} for a two-dimensional domain
discretized into triangles. The nodes of that mesh are the summits
of the triangles (note that, while common, this situation is not
the rule, and the middle of the edges, the center of gravity of the
triangles, $\ldots$ can also be nodes).
The effect of numbering is demonstrated in Figures \ref{ref_mesh}-a
and \ref{ref_mesh}-b, where
the same mesh can give a bandwidth of 18 or 26, depending
on the order in which the nodes are considered.\\

The only useful data for mesh numbering is, for each node,  the list
of all neighbors.
For instance, the mesh of Figure \ref{ref_mesh}-a can hence be viewed
as

\noindent
{\tt (2 3 4)(1 3 5)(1 2 4 5 6)(1 3 6 7)(2 3 6 8)(3 4 5 7 8)(4 6 8)(5 6 7)}

Once an order for the nodes is chosen, 
the profile of the matrix can be constructed: an upper bound
\footnote{Some of these candidates to be non-zero can
actually be null, depending on the actual formulation
of the PDEs. But this rare eventuality will not be considered
in the general framework developed here.} 
of the fitness function is given by the following equation (\ref{bw_id}),
where $N$ is the number of nodes, and, for each $i$ in $[1,N]$,
${\cal N}(i)$ is the set of neighbors of node $i$:

\begin{equation}
\label{bw_id}
 {\cal B} = \sum_{i=1}^{i=N} 
\begin{array}[t]{cc}
{min}\\
^{(j > i) ; (j \in {\cal N}(i))}
\end{array}
 (j-i) 
\end{equation}

Note that due to symmetry, only the contribution to the bandwidth of 
the upper part of the matrix is considered.

For the simple case of the mesh of figure \ref{ref_mesh},
the bandwidth is $3+3+3+3+3+2+1=18$ for the order (a) and $7+6+5+4+2+1+1=26$
for the order (b).

The goal of mesh renumbering is to find the order of the nodes
(i.e. a permutation on $[1,N]$) that minimizes the bandwidth $\cal B$.

\subsection{State of the art}
\label{Gibbs}
The problem of mesh numbering is clearly a $NP-complete$ 
combinatorial problem, as the search space is the space
of permutations in $[1,n]$. Hence, no exact deterministic 
method can be hoped for.
Numerical scientists have paid much attention to that
problem, developing heuristic methods.
The favorite method nowadays, used in most FEM software
packages, is the so-called {\em Gibbs method} presented in detail in 
\cite{George}).
It performs three successive steps, and use the graph $G$ representing
the ``neighbor to'' relationship between nodes.
\begin{itemize}
\item Find both ends of the numbering. They are chosen such that,
first, their distance (in term of minimal path in $G$ 
joining them) is maximal, and second, 
each one has the lowest possible number of neighbors.
\item Make a partition of graph $G$ into layers containing
nodes at the same distance from the origins.
\item Number the nodes breadth-first, starting from one origin (i.e. numbering 
nodes layer by layer).
\end{itemize}

As Gibbs method is used in all FEM packages, the minimum requirement
for any other algorithm to be considered interesting by numerical
scientists is to obtain better results than those of Gibbs.
Hence, in the following, all results of performance will be given 
relatively to those of Gibbs.

\section{The evolutionary algorithm}
\label{EA}

\subsection{MNP is not a TSP}
Combinatorial optimization is a domain where evolutionary algorithms
have encountered some successes, when compared to state-of-the-art
heuristic methods. Probably the most studied combinatorial
optimization problem is the Traveling Salesman Problem (TSP). In both
cases (MNP and TSP), the search space is that of all permutations of
$N$ elements.

However, a first obvious difference is
that both the starting point and the ``direction'' 
of numbering are discriminant in the MNP while they are not in the
TSP. As a consequence, the size of the search space in the MNP
is $n!$ while it is $(n-1)!/2$ for the TSP. On the other hand, 
no {\em degeneracy} (see \cite{Radcliffe-variance}) is present in the
{\em permutation} representation  for the MNP (see below).

Second, the MNP is not so easily {\em decomposable}. 
A whole part of a solution of the TSP can be modified without modifying
the remaining of the tour. In the MNP on the other hand, 
the propagation of any modification has consequences
on all geometrical neighbors of all modified nodes.

As a consequence, the useful notion of ``neighbor''
is totally different from one problem to another:
In the TSP, two towns are usually called neighbors if they are
visited one after the other on the tour at hand; In the 
MNP, the ``neighbor of'' relationship is absolute geometrical data,
independent of any order. 

Practically, a major difference between both problems is the absence
of any local well known optimization algorithm for the MNP: as quoted in
the introduction, the
best results obtained so far by evolutionary algorithms on the TSP are
due to hybrid ``memetic''
algorithms searching the space of local optima with respect to 
 a local optimization procedures
(e.g. the 2-opt or 4-opt procedures in
\cite{TSP:PPSN96,Merz:ICEC97,Asparagos:ICEC97}). Such strategy cannot
be reproduced for the MNP.

\subsection{Representation}
\label{representation}
One important consequence of these differences is that the concepts of
``edges'' or ``corners'' \cite{Radcliffe-variance}, known to be of
utter importance for the TSP, do not  
play any role in the MNP.

Two representations will be experimented with in that paper:
\begin{itemize}
\item In the {\em permutation} representation,
a permutation is represented by the sequential list of the numbers of all
nodes. Note that this representation 
relies on a predefined order of the nodes (corresponding to
the identity permutation $(1,2,3,\ldots,N)$). Some 
consequences are discussed in section \ref{initialization}.
\item All permutations can be decomposed in a sequence of 
transpositions. Moreover, all sequences of transpositions
make a unique valid permutation. The {\em transposition}
representation describes a permutation as an ordered list
of transpositions. Note that this representation is
highly degenerate: many genotypes correspond to the same phenotype.
\end{itemize}

\subsection{Crossover operators}

Four crossover operators have been tested, from general-purpose
operators to MNP-specific crossovers.
\begin{itemize}
\item The transposition crossover is a straightforward crossover
for the transposition representation: It exchanges portions of the
transposition lists of both parents. As any combination of transpositions
make a valid permutation, it directly generates valid permutations.
\item The {\em edge crossover} is a general-purpose operator
for permutation problems, designed and tested on the TSP
problem \cite{Whitley_TSP89}. 
No specific knowledge about the problem is used, but the underlying
assumption is that edges are the important features in the
permutations. All experiments using the edge crossover gave lousy
results, thus confirming that edges do not play for the MNP the
important role they have in the TSP. That crossover will not be
mentioned any more here. 
\item The {\em breadth-first} crossover is based on the heuristic
technique described in section \ref{Gibbs}, but uses additional information
from both parents to generate one offspring.

A starting point is randomly chosen, is given the number it has in
parent $A$, and becomes the current node.
All neighbors of the current node are numbered in turn, being awarded 
an unoccupied number, before being put in a FIFO stack. 
To number the neighbor $N$ of a node $M$ already numbered $i_M$,
the differences $\Delta_A$ and $\Delta_B$ of the numbers of nodes
$M$ and $N$ respectively in parent $A$ and parent $B$ are computed.
If the number $i_M + min(\Delta_A,\Delta_B)$ is free, it is given to
node $N$. Otherwise, if the number $i_M + Max(\Delta_A,\Delta_B)$
is free, it is given to node $N$. Otherwise, the number closest 
to $i_M$ which is not yet used is given to node $N$.

All nodes are processed once, and are given a yet-unattributed
number, generating a valid permutation. This crossover tries
as much as possible to reproduce around each node the best local 
numbering among those of both parents.

\item The {\em difference} crossover was designed to both preserve the
diversity in the population and try to locally minimize the fitness
function. 
From both permutation representations of the parents, all nodes having
the same number are given that common number in the offspring
permutation. Further, all remaining nodes in turn (in a random order)
are given the number which is the closest possible from all
its neighbors numbers. 
\end{itemize}

\subsection{Mutation operators}
\label{mutation}
Here again, both general-purpose and problem specific
mutation operators were tested.

\begin{itemize}
\item The minimal mutation for the permutation representation
is the exchange of numbers between two randomly chosen nodes
(i.e. the application of a transposition), thereafter termed
{\em transposition mutation}. 
Its strength can be adjusted by repeated application 
during the same mutation operation.
\item The {\em neighbor transposition mutation} is a 
slight modification of the above operator: after the first node
has been chosen randomly, it is exchanged with one of its neighbors.
\item The {\em neighbor permutation mutation} is a further step
in the direction of using neighbor information to perform mutation:
a node is randomly chosen, and a random permutation among the numbers
of all its neighbors is performed.
\item The {\em inversion mutation} reverts some part of the permutation.
A special case of the inversion mutation is when the whole permutation
is inverted. This operator was found useful in Gibbs method (section
\ref{Gibbs}). 
\item The choice of the origin of the numbering is important in the MNP.
The {\em origin mutation} was designed to handle this issue.
An integer $i$ is randomly chosen, and all numbers $j$ are replaced by
either $i+j \mbox{ mod } n$ either $n-(i+j) \mbox{ mod } n$, on a local minimization
argument (the brute translation only gives birth to usually
very bad numbering when the old origins met).
\end{itemize}

\subsection{Initialization procedure}
\label{initialization}
The standard way to initialize a population
is to perform a uniform sampling of the genotype space.
However, alternative specific ways
have been proposed: for instance, the greedy heuristic for the TSP
(from a random initial town, chose the nearest town not yet visited)
constructs individuals with a tour length of about 20\% 
more than the optimal value on average \cite{TSP:PPSN96}.
In the same line, three different initialization procedures have been
tested for the MNP. 

\begin{itemize}
\item The standard uniform sampling of the standard representation
(section \ref{representation}) was the first obvious choice. It is termed
{\em random initialization}. 
\item As the state-of-the-art solution is given by the output of the
Gibbs method (section \ref{Gibbs}), the {\em Gibbs initialization}
performs only slight perturbations of the Gibbs solution to generate
the initial permutations. This is achieved using the transposition
representation (section \ref{representation}) taking the Gibbs
order as reference, and allowing only a small number of
transpositions. Note that in this case, the original Gibbs result is
included as the first individual of the population.
\item The trouble with the above Gibbs initialization is a very strong
bias toward solution quite similar to Gibbs. If higher optima
are located in very different regions of the permutation space,
they will probably not be found. To address that issue, the {\em point
initialization} was designed, based on some breadth-first
heuristic similar to 
Gibbs', but using a random starting point (Gibbs process is very sensitive
to the choice of the initial point). 
\end{itemize}

The average bandwidth of permutation drawn using the  uniform initialization
is of course quite large (more than 6 times that of Gibbs method). 
The point initialization gives much better individuals: their
bandwidth range from 20\% to 100\% above Gibbs results while Gibbs
initialization stays between 0\% and 15\% above Gibbs order.

\section{Experimental Results}
\label{results}
\subsection{The meshes}
Three meshes have been used to test the algorithm described above: the
{\em small} mesh has 164 nodes, the {\em medium}
one has 1544 nodes (see Figure \ref{test-mesh}) and the {\em large}
mesh has 5453 nodes.
The computational cost of one fitness function evaluation increases
linearly with the size of the mesh. Hence, most initial experiments
were performed in the {\em small} mesh. The validation of clear
tendencies was then carried on for confirmation on the {\em medium}
mesh. The best combination of parameters were finally tested on the
{\em large} mesh, as its size is becoming to be  of some interest for
real world application.

\subsection{Experimental settings}

Two basic evolution schemes were used: a standard GA, with linear ranking and
elitist generational replacement; a ($\mu$+$\lambda$)-ES scheme, in
which all $\mu$ parents give birth to $\lambda$ offspring, the best
$\mu$ of the $\mu + \lambda$ parents + offspring becoming the parents
of the next generation. 
The first series of experiments were performed on the
{\small} mesh, using a population sizes of 50
for the GA scheme, and a (7+50)-ES.
Both
schemes were tested with the same combinations of $P_c$-$P_m$,
(crossover rate - mutation rate).
If an individual undergoes  crossover, 
a mate is selected and only one
offspring is generated. One single type of crossover was made possible
at each run.
The resulting offspring (or the initial individual if
no crossover occurred) then undergoes mutation with probability $P_m$.
If mutation happens, one mutation operator is chosen according to 
user-defined weights, and, in the case of transposition mutations, one
single transposition is performed.

\subsection{First tuning on the small mesh}
As said above, intensive experiments were performed on the small mesh.
During those experiments, 
all values for $P_c$ and $P_m$ between 0 and 1 by 0.1 steps were tried
independently, for all possible crossover operators. Those runs were
allowed 150000 fitness evaluations (unless otherwise mentioned).\\

Preliminary runs were performed to tune the mutation weights
(\cite{Sourd:Option96}).
The weights for the 5 mutations
were first set equal. Then  a close look at the
types of mutations that the best individual in the population was
submitted to, along different runs, allowed to eliminate both the {\em
inversion mutation} and the {\em origin mutation}. 
Only the  {\em random transposition mutation}, the {\em neighbor
transposition mutation} and the {\em neighbor permutation mutation}
proved useful, and their weights $PmRand$, $PmNeighbor$ and $PmAround$  were
set to values between 0 and 1, their sum being equal to
1. 

\subsubsection{Initialization procedures}
As could be predicted, the best on-line results were obtained for the
{\em Gibb's initialization}, as its starting point was rather better than
both other. But whereas the {\em point initialization} almost caught
up (as will be seen in forthcoming section \ref{res-mutations}), the
{\em random initialization} stayed far beyond, even when allowed
ten times the number of fitness evaluations. So only the {\em Gibb's}
and the {\em point} initializations will be considered in the following.

\begin{center}
\begin{figure}
\begin{tabular}{ccc}
\psfig{figure=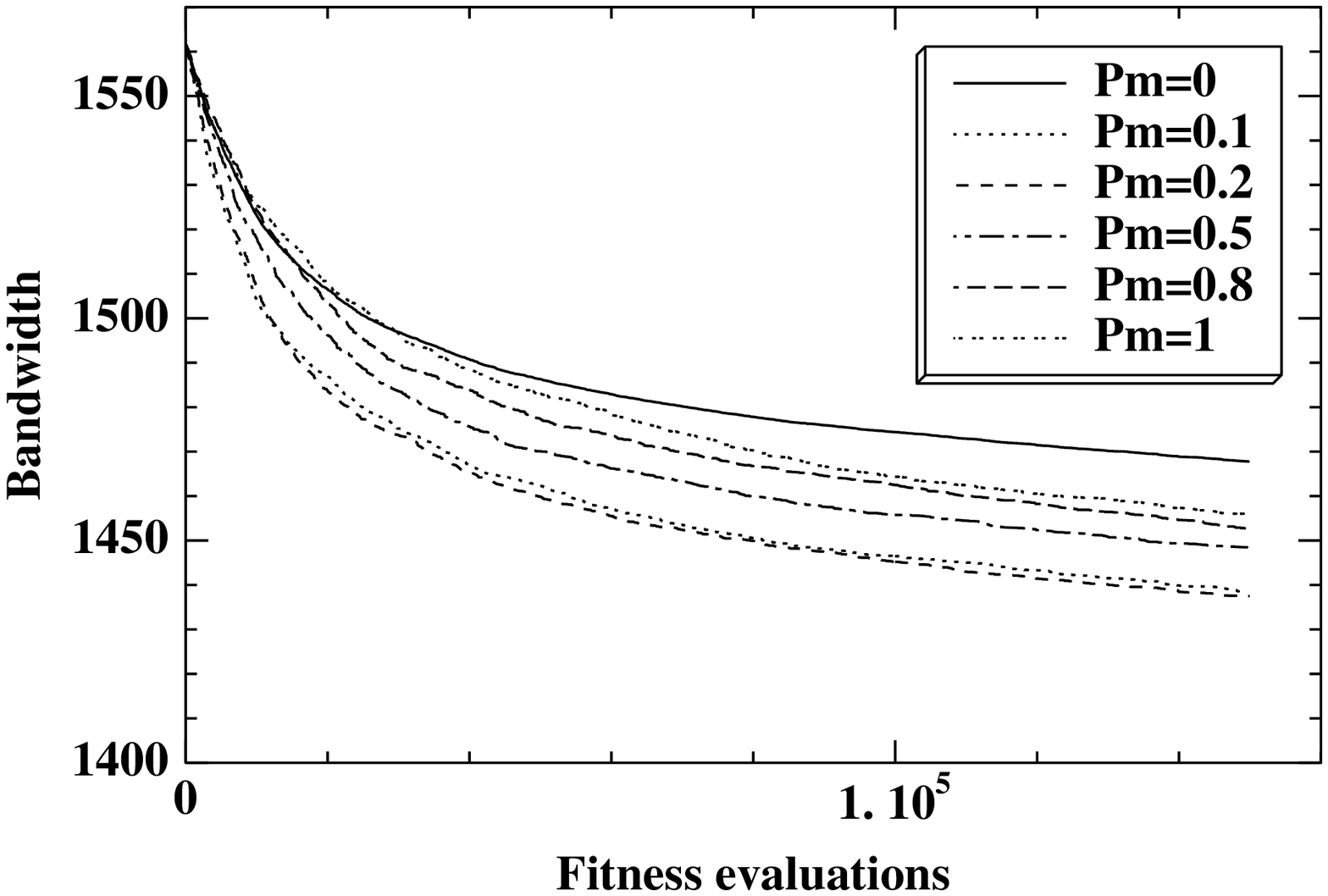,width=5cm} & ~~~~~&
\psfig{figure=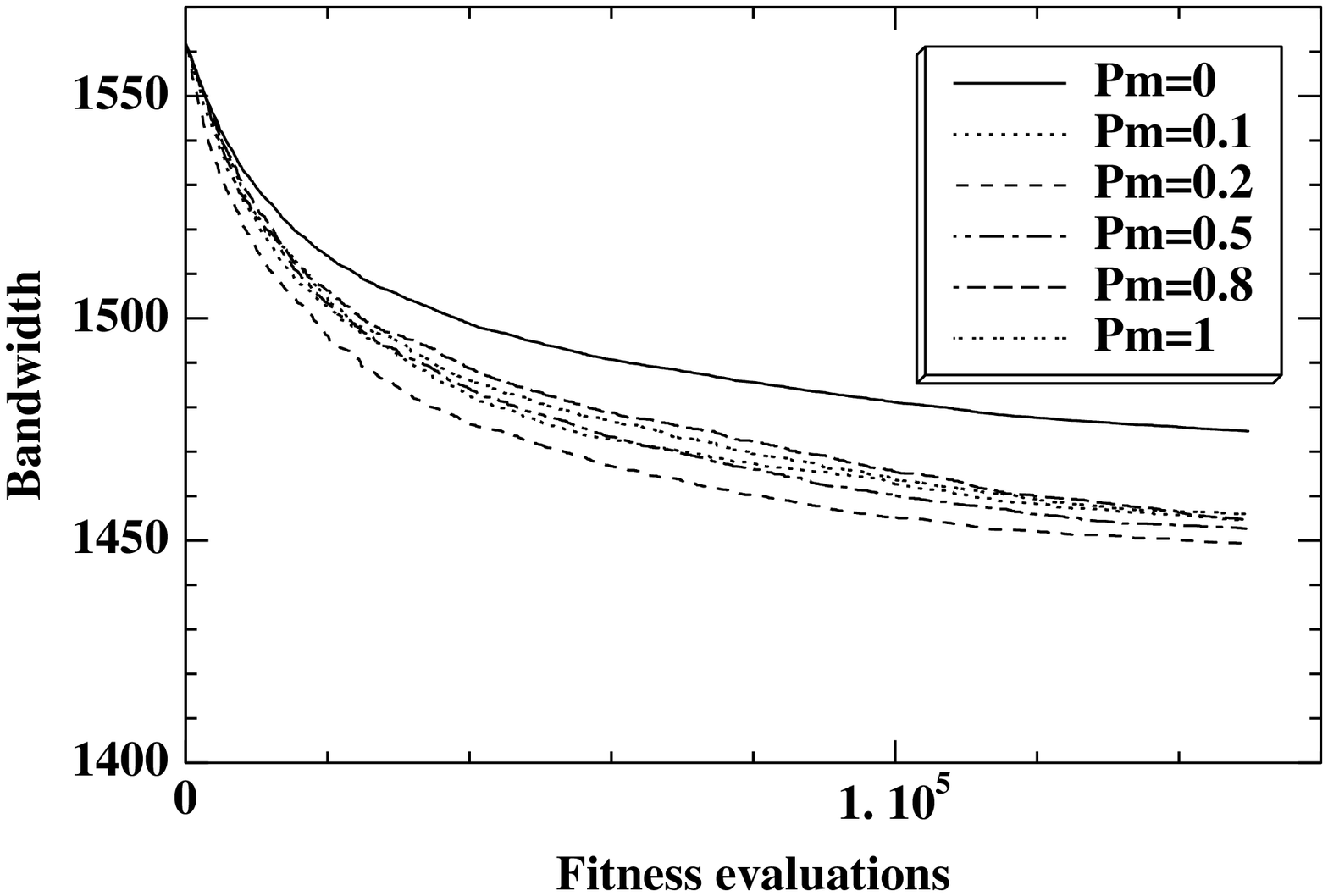,width=5cm} \\
(a) {\em GA with breadth-first crossover} & &
(b) {\em GA with difference crossover} \\
~\\
\psfig{figure=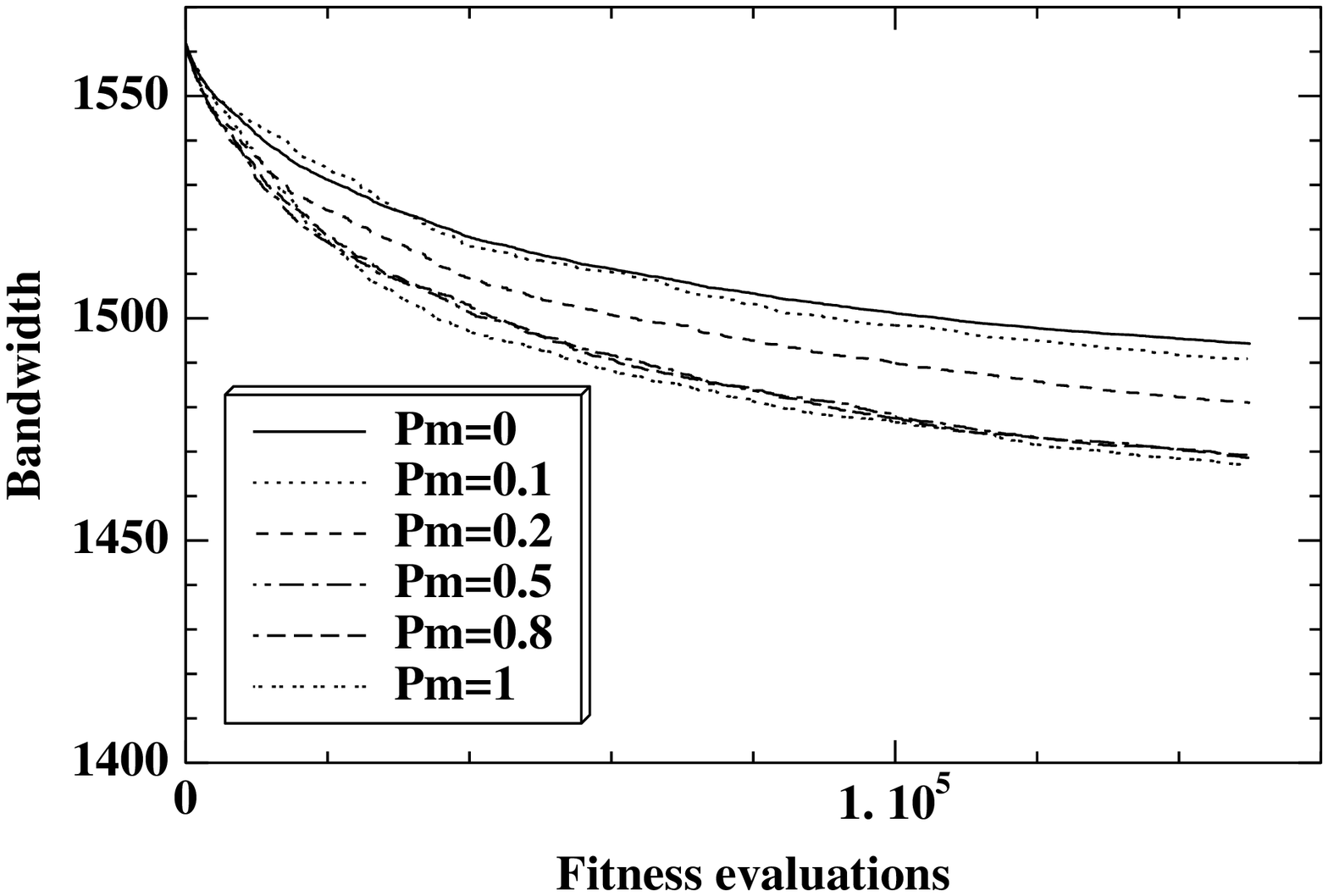,width=5cm} & ~~~~~&
\psfig{figure=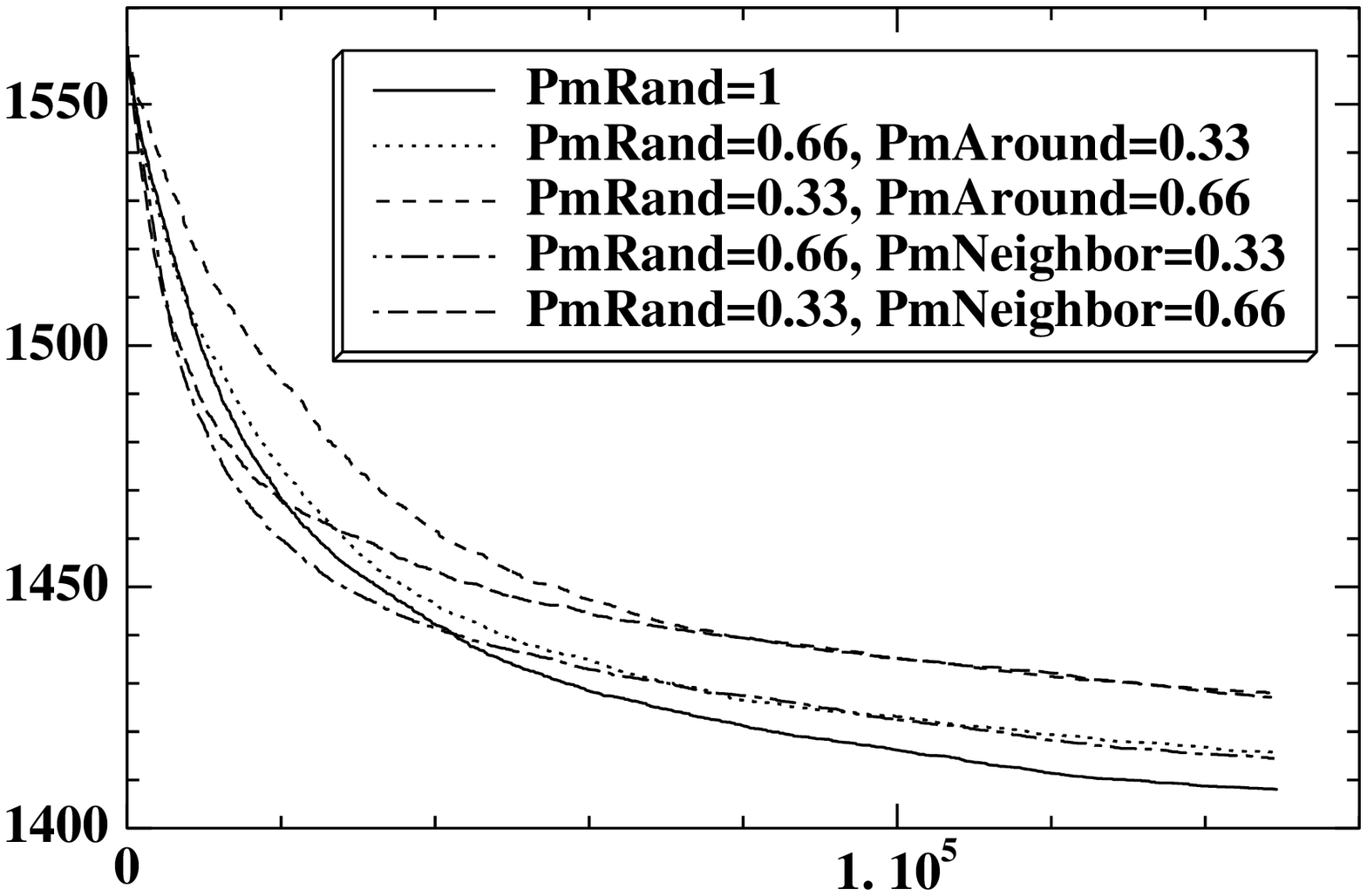,width=5cm} \\
(c) {\em GA with transposition crossover} & &
(d) {\em (7+50)-ES - no crossover}\\
\end{tabular}
\caption{\em Typical on-line results (averages over 21 runs) on the
small mesh: for the three GA runs, 
population size is 50, $P_c$ is 0.6 and $P_m$ is as indicated;
for the (7+50)-ES runs, $P_m$=1; in all cases, whenever an individual
is mutated, PmRand=0.66 and PmNeighbor=0.33.}
\label{res-small}
\end{figure}
\end{center}

\vspace{-1.5cm}

\subsubsection{Crossover operators and evolution schemes}
The first experiments aimed at comparing the crossover operators, and
adjusting the crossover rate $P_c$. 
A common feature could be observed for all three crossovers: 
when using the GA scheme, and except for high values of
$P_m$, for which it made no 
significant difference, the results decreased when $P_c$ decreased
from 1 to 0.6 or 0.5. Moreover, and almost independently of the
settings of the mutation weights, the best results were obtained with
the (7+50)-ES scheme, with $P_c=0$ (higher values of $P_c$ for the ES
scheme performed rather poorly - but this might be because of the
rather small population size of 7).

When it comes to compare the crossover operators, the results of
Figure \ref{res-small} (a), (b) and (c), are what could be expected:
the ``blind'' {\em
transposition crossover} (c) performs rather poorly -- and gets its best
results for the highest values of $P_m$ which seems to indicate that
it is really not helping much. On the opposite, both other
operators, that do incorporate some domain knowledge,  get their
highest performances for $P_m=0.2$ and $P_m=0.1$. 
Moreover, the 
{\em breadth-first crossover} (a) performs better than the {\em
difference crossover} (b) (and the difference is statistically
significant with 99\% confidence T-test after 300000 fitness computations). 

A last argument favoring the abandon of crossover operators is the
extra cost they require, as based on local optimization heuristics of
complexity $o(N)$.
For instance the total CPU time is increased by a factor around $4$
between runs with $P_c=0$ and $P_c=1$ (from 3 to 13mn on average for
300000 evaluation runs on a Pentium P200).

\begin{center}
\begin{figure}
\begin{tabular}{ccc}
\psfig{figure=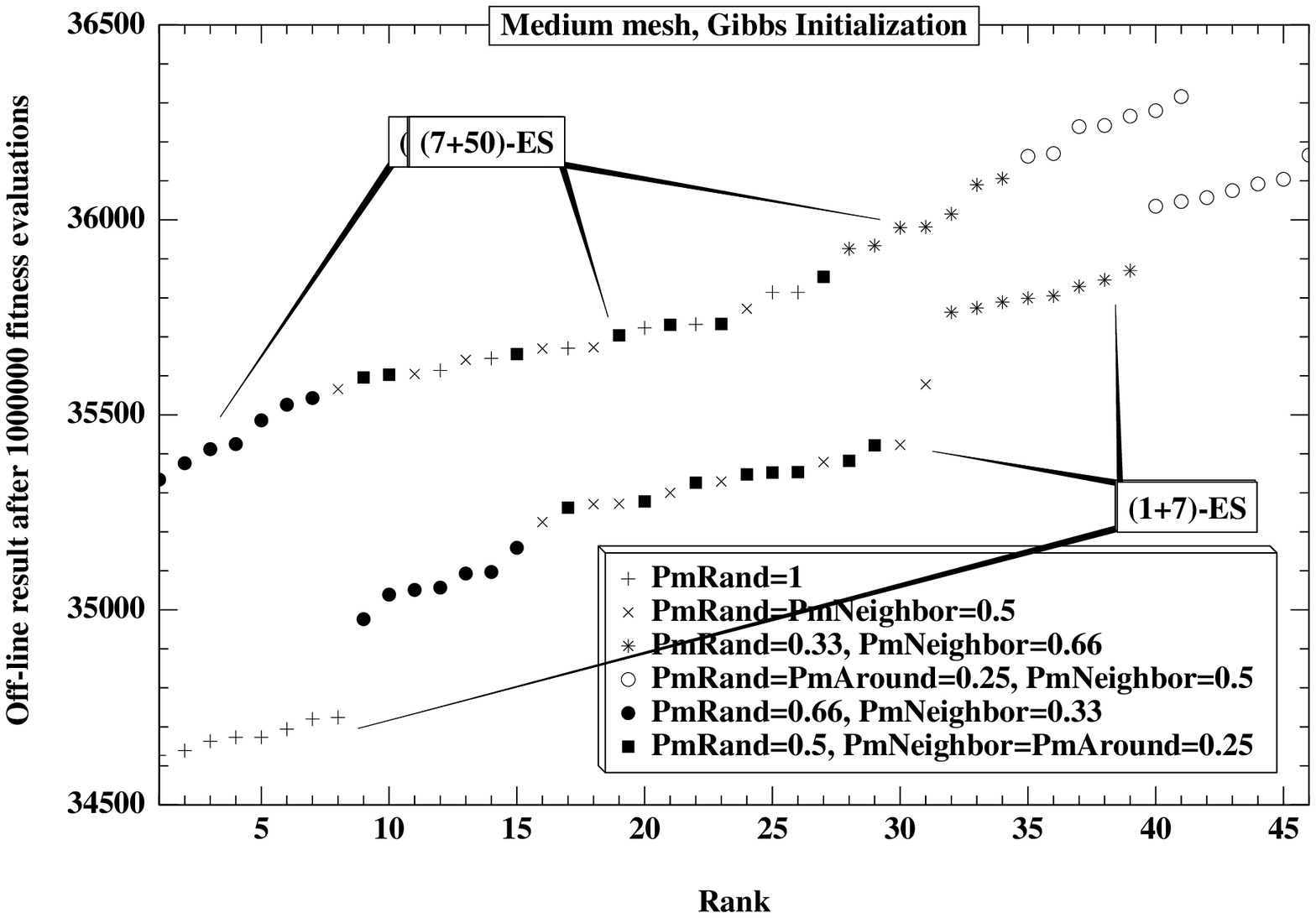,width=6cm} & ~~~~~&
\psfig{figure=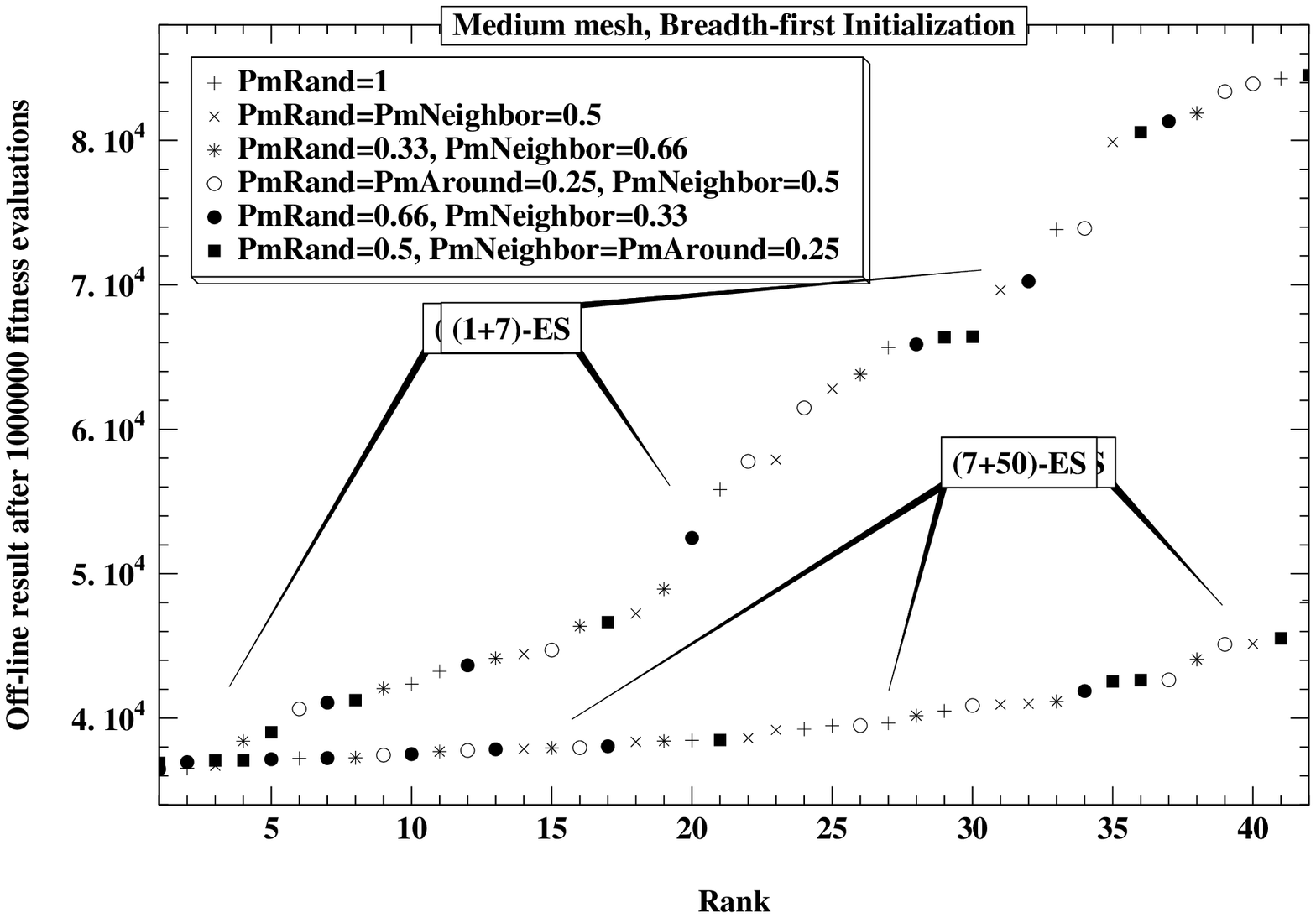,width=6cm} \\
(a) {\em Gibbs initialization} & ~&
(b) {\em Breadth-first initialization} 
\end{tabular}
\caption{\em Off-line results of both (7+50)-ES and (1+7)-ES
after 1000000 fitness evaluations, for different settings of the
mutation weights.}
\label{res-medium}
\end{figure}
\end{center}

\vspace{-1.5cm}

\subsection{Mutations and population size}
\label{res-mutations}
This section presents some results obtained on the {\em medium} and
{\em large} meshes, using mutation operators only inside some $ES$
evolution scheme. At first, the goal of these experiments was to sort
out the usefulness of problem-specific knowledge, in the initialization
procedure and in the mutations operators. 
Bur it rapidly turned out that the population size also was a very
important parameter. 

Figure \ref{res-medium} witnesses the surprising
results obtained on the medium mesh: each dot indicates the best
fitness reached after 1 million fitness evaluations of a single run of
the evolutionary algorithm. The different shapes of the dots represent
different settings of the relative mutation rates $PmRand$,
$PmNeighbor$ and $PmAround$. Figure \ref{res-medium}-a shows the runs
that used the Gibbs initialization procedure while Figure
\ref{res-medium}-b those who used the breadth-first initialization
procedure. On each Figure, the two distinct sets of points correspond
to the $(1+7)$- and the $(7+50)$-ES schemes, as indicated. Note that
all trials with larger population sizes were unsuccessful.

Some clear conclusions can be drawn from these results.
First, the overall best results are obtained for the $(1+7)$-ES scheme
starting from the Gibbs initialization and using the ``blind'' {\em
transposition mutation} only (Figure \ref{res-medium}-a) (see section
\ref{discussion} for a detailed comparison with the results of Gibbs
method). But comparing the results between Gibbs and breadth-first
initialization procedures on the one hand, and $(1+7)$- and $(7+50)$-ES
schemes on the other hand gave some unexpected results:

\begin{itemize}
\item Whereas the $(1+7)$-ES scheme
consistently outperformed the the $(7+50)$-ES scheme when using 
Gibbs initialization, the reverse is true for the breadth-first
initialization;
\item Whereas different settings of the mutation rates gave rather
different results for the Gibbs runs, this does not seem to be so
clear for the breadth-first runs;
\item whatever the initialization, the $(7+50)$ results are more
stable than the $(1+7)$ results with respect to mutation rates; This
is striking on the breadth-first plot, but also true on the other
plot.
\item the worst results are obtained for the lowest values of $PmRand$
(the results with even lower values of $PmRand$ are not presented
here, but were very poor). However, domain
specific mutation operators (see the black circles and squares,
compared to the ``+'' dots) are
more efficient with the $(7+50)$ 
scheme than with the $(1+7)$ scheme. This is specially visible on the
Gibbs runs.
\end{itemize}

Note that these tendancies were confirmed when the runs were allowed
more fitness evaluations (e.g. 10 millions, see forthcoming section
\ref{discussion}). 

Some tentative explanations can however be proposed, after noticing
that the $(1+7)$ scheme can be viewed more like a {\em depth-first}
algorithm while the $(7+50)$ scheme searches in a more {\em
breadth-first} manner. 

So, using Gibbs initialization probably tights the population in a
very limited area of the search space from which it is useless to try to
escape: this favors the performance of depth-first search, as breadth
oriented search does not have the possibility to jump to other
promising regions. Moreover, it seems that the transpositions of
neighbor nodes in a depth-first search does not 
allow large enough moves to easily escape local optima, resulting in
the best results for the pure random mutation for the $(1+7)$. 

On the other hand, successive local moves have greater chances
of survival in the $(7+50)$ breadth search, and give the best results
in that case (though some random mutation are still needed). And when
it comes to a more widely spread population after the breadth-first
initialization, the breadth-first search demonstrates better results
by being able to use more efficiently in that case the domain
neighboring information.\\

This situation suggests further directions of research: First, the
$(1+7)$-ES scheme with pure random mutation resembles some sort of
Tabu search \cite{Glover:77,Glover:95}, and so might be greatly improved
by adding some memory to it, either deterministically, like in
standard Tabu search, or stochastically, as proposed in
\cite{Nous:ICGA97,Peyral:EA97}.

Second, more than one neighbor transposition  seems necessary to
generate improvements. Hence, the number of transpositions should not
be forced to 1, and can be made
either self-adaptive, with the same problems than in the
case of integer variables in Evolution Strategies \cite{Baeck:EP95},
or exogenously adaptive, as proposed in \cite{Baeck-hyperbole} where
some hyperbolic decreasing law is used for the mutation rate along
generations.

\subsection{Evolutionary mesh numbering vs Gibbs method}
\label{discussion}
But apart from optimizing the evolutionary algorithm itself on the
MNP, a critical point is whether evolutionary mesh numbering can
compete with the Gibbs method. Of course, if the Gibbs initialization
is used, as the original Gibbs numbering is included in the initial
population, any improvement is in fact giving a better result than the
Gibbs method. But how interesting is that improvement, especially when
the computational cost is taken into account?\\

On the medium mesh (1545 nodes), the bandwidth using Gibbs method is
39055. As can be seen on Figure \ref{res-medium}, the best result of
the $(1+7)$-ES with pure random mutation after 1 millions evaluations
is 34604, i.e. an improvement of 11.4\%. 
The computational cost of one run is 15-20mn (on a Pentium 200Mhz
Linux workstation), to be compared to the
20s seconds of Gibbs method!

If the maximum number of  function evaluations is set to 10 millions,
the best result in the same conditions is 32905 (i.e. 15.75\%
improvement), with an average over 21 runs of 33152 (15.11\%).
Of course, in that latter case, the computational cost is 10 times
larger $\ldots$\\

The results on the large mesh (5453 nodes) follow the same lines,
though only the combinations of parameters found optimal on the medium 
mesh were experimented with, due to the computational cost (around
12-14 hours for 10 Millions evaluations). 

From a Gibbs bandwidth of 287925 (obtained in about one minute), 
the best result for the $(1+7)$-ES  with only random mutation was
257623 (10.52\%) 
while the average of 6 runs was 258113 (10.35\%).
On the other hand, the best for $(7+50)$-ES was 262800 (8.73\%), 
the average being 263963.43 (8,32\%), 
with $PmRand = 0.66$ and  $PmVois = 0.33$ (best parameters of Figure
\ref{res-medium}-a). \\

The first a priori conclusion is that the computational cost of
the evolutionary algorithm makes it useless in practical situations: an
improvement of between 
10 and 15\% requires a computational power of several order of
magnitude larger.
This quick conclusion must however be moderated: First, due to the
quadratic dependency of the computational cost of the matrix inversion
in term of the bandwidth, the actual gain in computing time is around
35\% for a 15\% bandwidth decrease. 
And second, many meshes used nowadays in industry require a few
months of manpower to be built and validated. So 24 more hours of
computation for a better numbering is relatively low increase in
cost. And if the mesh is then put in an exploitation environment, and
is used in  several thousands of different Finite Element Analyses,
then the overall gain might be in favor of getting a really better
numbering, even at what a priori seems a high computational cost.
But of course this means that meshes of up to a few thousands nodes
can be handled by evolutionary algorithms.

It is nevertheless important to notice that the computation of the
bandwidth can be greatly optimized. The present algorithm was designed
to be very general, handling any possible operator. Hence it always computes
the fitness from scratch. However, in the case where only a small
number of transpositions are performed, the variation of the fitness
could be computed, by examining only the neighbors of the transposed nodes.

\section{Conclusion}
\label{conclusion}

We have presented feasibility results for the application of
Evolutionary Computation to the  problem of mesh
numbering. Our best results outperform the state-of-the-art Gibbs
method by 10 to 15\% on the two test meshes used (with 1545 and 5453
nodes respectively). Whereas these sizes would appear fairly high for
TSP problems for instance, they are still small figures with respect
to real-world problems, where hundreds of thousands of nodes are
frequent.

From the Evolutionary point of view, two issues should be highlighted.
First, though both general-purpose and domain-specific crossover
operators were tried, none proved efficient. A possible further trial
could be to use more global geometrical information (e.g. divide the
mesh into some connected components, and exchange the relative orders
of such blocks, in the line of \cite{KS_EA95}).

Second, the overall best results were obtained by a $(1+7)$-ES using
pure random transposition mutation and starting from an initial
population made of slightly perturbed Gibbs meshes. This which might
be an indication that 
other heuristic local search methods (e.g. Tabu search) might be
better suited to the MNP. However, as discussed in section
\ref{res-medium}, some hints make us believe that there is still a
large room for improvement using evolutionary ideas: on the one hand,
the problem-specific mutations proved useful for the $(7+50)$-ES,
indicating that we might not have make good usage of the domain
knowledge; on the other hand, the $(7+50)$-ES (with problem-specific
mutation) outperformed the $(1+7)$-ES when the initial population was
not limited to modified Gibbs meshes: our hope is that
 starting from totally different parts of the search space could provide
much better results in some particular situations ... which still remain to
be identified. 
But in those yet hypothetical  cases, the relevance of the
evolutionary approach for the MNP would be clear.

\small

\end{document}